\documentclass[10pt]{article}
\pdfoutput=1
\usepackage{amsmath}
\usepackage{amssymb,amsmath,amsthm}
\usepackage{graphicx}
\usepackage{amssymb}
\usepackage{cite}
\usepackage{cite}
\usepackage{xfrac}
\usepackage{url}
\usepackage{url}

\usepackage{lipsum}
\newcommand\blfootnote[1]{%
  \begingroup
  \renewcommand\thefootnote{}\footnote{#1}%
  \addtocounter{footnote}{-1}%
  \endgroup
}

\newtheorem{thm}{Theorem}
\newtheorem{lem}{Lemma}
\newtheorem{prop}{Proposition}

\newtheorem{defin}{Definition}
\newenvironment{example}[2][Example]{\begin{trivlist}
\item[\hskip \labelsep {\bfseries #1}\hskip \labelsep {\bfseries #2}]}{\end{trivlist}}
\newenvironment{rem}[1][Remark]{\begin{trivlist}
\item[\hskip \labelsep {\bfseries #1}]}{\end{trivlist}}
\newcommand{\bd}{\begin{defin}}
\newcommand{\ed}{\end{defin}}
\newcommand{\bl}{\begin{lem}}
\newcommand{\el}{\end{lem}}
\newcommand{\br}{\begin{rem}}
\newcommand{\er}{\end{rem}}
\newcommand{\bt}{\begin{thm}}
\newcommand{\et}{\end{thm}}
\newcommand{\be}{\begin{example}}
\newcommand{\ee}{\end{example}}
\newcommand{\bp}{\begin{prop}}
\newcommand{\ep}{\end{prop}}

\newcommand{\bqn}{\begin{eqnarray}}
\newcommand{\eqn}{\end{eqnarray}}

\def\1{\hbox{\rm\rlap {1}\hskip.03in{\rom I}}}
\def\Bbbone{{\rm1\mathchoice{\kern-0.25em}{\kern-0.25em}
{\kern-0.2em}{\kern-0.2em}I}}

\begin{document}
\title{{\LARGE{Persistent Homology of Delay Embeddings and its Application to Wheeze Detection \\ }}
}
\author{   \Large{{Saba~Emrani, Thanos~Gentimis and Hamid~Krim}}  \\
}
\maketitle
\begin{abstract}
We propose a new approach to detect and quantify the periodic structure of dynamical systems using topological methods. We propose to use delay-coordinate embedding as a tool to detect the presence of harmonic structures
by using persistent homology for robust analysis of point clouds of delay-coordinate embeddings. To discover the proper delay, we propose
an autocorrelation like (ACL) function of the signals, and apply the introduced topological approach to analyze breathing sound signals for wheeze detection.
Experiments have been carried out to substantiate the capabilities of the proposed method.
\end{abstract}
\section{Introduction}
\blfootnote{The authors are with the Department of Electrical and Computer Engineering, North Carolina State University, Raleigh, NC 27606 USA (emails: semrani, agentim, ahk@ncsu.edu).\\
The work is supported by NSF, EEC-1160483 CCF-1217874.} 
As previously noted in \cite{Carlsson08}, persistent homology can be a powerful topological  approach for analyzing large data sets. Topological tools, such as barcodes and Betti numbers of the spaces can be employed on data samples from dynamical systems to extract various characteristics. Time-delay coordinate embedding has mostly been used in the analysis of dynamical systems \cite{Takens}. 
Takens has proved an important theorem which states that almost every time-delay embedding of a time series can recover the underlying dynamics of a system \cite{Takens}. There exists a large body of literature on the application of the delay-coordinate embedding technique to  chaotic attractors  \cite{Abar, Kantz, Han11}. 
Knudson et al. \cite{Kevin} used three-dimensional time delay embeddings of human speech signals by approximating the data sets by simplicial complexes and by analyzing their persistent homology.
In this study, we propose the use of topological methods to discover almost periodic structures in the delay-coordinate embeddings of data samples of dynamical systems with little prior associated information. 
An autocorrelation like (ACL) function is used to select an appropriate value for delay. This novel approach for signals and systems analysis is general, and is robust to missing data points since we describe the global topological structure of the point cloud rather than their local geometric behavior. Moreover, our proposed method is computationally fast, since only a small number of data points are extracted from the whole point cloud, and are used.

Wheezes are abnormal lung sounds, and usually imply obstructive airway diseases. The most important characteristic of wheeze signals, which is the key component of this study, is their harmonic behavior. Wheeze detection problem is studied in \cite{Taplidou10, Homs04, Taplidou07}. We will apply our proposed approach for faster and more accurate wheeze detection, which is also robust to missing data points.
We propose a model for wheeze signals as a piecewise sinusoidal function with time varying frequencies and amplitudes. We show that wheeze signals are at $\epsilon$ distance of their corresponding introduced model. Accordingly, we can use wheeze signals as a good illustration of the topological analysis of time-delay embedding.  

The remainder of the paper is organized as follows: 
Time delay embedding is introduced in its generality in Section II, where an autocorrelation-like function is also proposed as a way to select a time delay. The topological characterization of the delay embeddings using persistent homology is then analysed. A model for wheeze signals is proposed and validated in Section III.  In Section IV, experimental results are given, and finally Section V concludes the paper.  
\section{Proposed Framework}
\subsection{Time Delay Embeddings}
We propose to use delay-coordinate embedding as a tool to construct a point cloud from digital signals and extract their periodic behavior.
The mathematical foundation of the delay-coordinate embedding method which embeds a scalar time series into an $m$-dimensional space, was proposed by Packard et al. \cite{Pack} and Takens \cite{Takens}. For each time series $\lbrace x_i\rbrace , i=1,2,...$, a representation of the delay coordinate embedding can be described as the following vector quantity of $m$ components:
\begin{equation}\label{gde}
  X_i=\left(x_i, x_{i + j},x_{i + 2j},...,x_{i + (m-1)j}\right), \quad X_i \in \mathbb{R}^m
\end{equation}
where $j$ is the index delay and $m$ is the embedding dimension.
If the sampling time is $T_s$, then the delay time $\tau$ is connected to the index delay $j$ by the equality $\tau=j\cdot T_s$. 
In this study, we use Equation (\ref{gde}) with  $m=2$ to detect almost periodicity and refer to it as 2-dimensional delay embedding. We will leave the study of $m=3$ for estimation of different frequencies present in the data set to up-coming papers. 
%
Discrete time sound signals can be considered as a series expressing the amplitude of the wave in volts at each time instance,i.e. $x(t_i):=x_i,i=1,2,...,k$ where $t_i=i\cdot T_s$. 
\subsection{Selecting Time Delay}\label{delay}
The time delay needs to be carefully chosen in order for the delay-coordinate components to be independent.
If it is too small, the delay embedding is compressed along the identity line. On the other hand, for too large delays, the adjacent components of the delay embedding may become irrelevant \cite{Cas91}. 
%
We propose to examine an autocorrelation-like function (ACL) to choose a proper delay. 
The ACL function of a non-stationary digital signal $x(t_i)$ is calculated as follows:
\begin{equation}\label{ac}
R_{xx}(t_i)=\sum_{1\le l\le k} x(t_i)\cdot x(t_l)\
\end{equation}
Clearly, peaks in the ACL function mark delay times at which the signal is comparatively highly similar to itself. 
According to experimental results, the appropriate interval for choosing delay time to best obtain informative delay embedding of signal $x(t)$ is $t_{c1}<\tau<t_{c2}$, where $t_{c1}$ and $t_{c2}$ are the first and second critical points of the ACL function $R_{xx}(t_i)$.
Note that for time varying non-stationary signals, a strict autocorrelation function cannot be defined and used for establishing periodicity nor for detecting the frequencies. However, we use Equation (\ref{ac}) as an ACL function strictly for delay selection.

\subsection{Persistent Homology and Point Cloud Subsampling}\label{landmark selection}
Computational topology uses topological invariants to distinguish objects and classify topologically equivalent spaces. Specifically, the ``Betti" numbers are ranks of a special type of topological invariants, called homology groups. In fact, $\beta_0$ gives the number of connected components of a topological space, $\beta_1$ measures the number of 1 dimensional topological holes, and $\beta_2$ counts the number of 2 dimensional topological ``voids".
The persistent homology algorithm uses a nested sequence of ``simplicial complexes" which approximates a certain space, called a filtration, to extract topological information about it. It characterizes the topology of these spaces using a collection of persistence intervals, or persistence intervals which are encapsulated in a structure called the persistent barcode. As we proceed in our filtration sequence of the space, each $n$-dimensional persistence interval $[t_b,t_d]$ can be associated to an n-dimensional hole, which appears at time $t_b$ (birth) and ``closes" at $t_d$. (death) \cite{Carlsson08}. We are focusing mainly on the first persistent Betti number as a means of distinguishing between wheeze and non-wheeze signals. 

For large datasets, using the whole set of points results in excessively large complexes with a high computational cost. 
Random and maxmin algorithms \cite{sil} can be used as two main options for the point cloud subsampling. 
After the point cloud subsampling, we utilize the Rips complex construction \cite{Edel10}.


The maxmin procedure is generally more reliable than random selection, and yields more evenly-spaced subsamples; albeit requiring a higher computational complexity. For the point cloud representing the data described in this study, the maxmin method shows a smaller hole compared to the random method since it emphasizes extremal points.
\section{Application to Wheeze Detection}
In this section, we first introduce a model for wheeze signals. We will then verify the proposed model by showing that the Hausdorff distance \cite{amme} of the graph of the model from the graph of the real wheeze signal is small.
The Hausdorff distance is one of the most commonly used metrics in shape analysis for comparing compact subsets. We exploit it by invoking a theorem to show that the persistence diagrams of the model and the wheeze signal are close.
\subsection{Wheeze Modeling}
A monophonic wheeze consists of a single note or several notes starting and ending at different times \cite{Lung09}.
The frequency of a monophonic wheeze signal is approximately a piecewise constant function of time. Accordingly, we propose a model for wheeze signals in the time domain defined as a continuous piecewise sinusoidal function with different periods and phases and a time varying amplitude, represented as
%
\begin{equation}\label{w}
  w(t)=\sum_{i=1}^n g_i(t),
\end{equation}
where
\begin{equation}\label{g}
  g_i(t) = \left\lbrace  \begin{array}{ccc}
                               w_i(t) &  & t_{i-1}\leq t < t_i,   \\
                               0 &  & \text{otherwise,}
                             \end{array}
                             \right.
\end{equation}
and $w_i \text{s}, i=1,2,...,n$ are defined as,
\begin{equation}\label{wi}
  w_i(t)=A(t) \sin\left(\frac{2 \pi}{T_i} t+\phi_i \right),
\end{equation}
where $A(t)$ is a nonzero continuous amplitude function. Also, to satisfy the continuity of $w(t)$, the phases $\phi_i$ should conform to the following condition
  $\phi _i=\phi _{i-1}+2\pi t_{i-1} (\frac{1}{T_{i-1}}-\frac{1}{T_{i}})$.

\begin{figure}[t]
\centering
\includegraphics[width=0.9\linewidth]{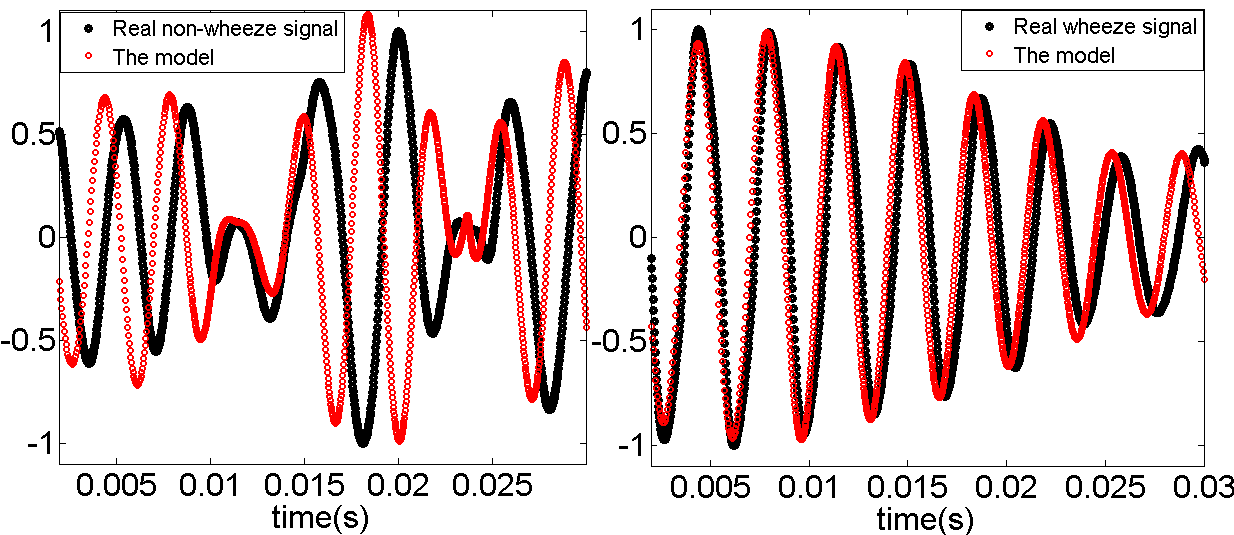} 
\caption{Comparison of the proposed model for a non-wheeze breathing sound signal (left) and a wheeze signal (right)} \label{mdl}
\end{figure}

To evaluate the effectiveness of the proposed model in representing wheeze signals, we show that the Hausdorff distance of the graph of the model from that of the recorded wheeze signal is sufficiently small.
Suppose that $I=\lbrace t_j : 1 \leq j \leq n, t_j \in [a,b], t_j=j \cdot T_s \rbrace$. Then for each discrete wheeze signal denoted by $s(t)$, which is the windowed version of a digital recorded wheeze sound in $[a,b]$, we can construct a function $w(t)$ defined using (\ref{w})-(\ref{wi}). If $A=\lbrace \left(t_j,s(t_j)\right): t_j \in I\rbrace \subseteq \mathbb{R}^2 $ and $B=\lbrace \left(t_j,w(t_j)\right): t_j \in I\rbrace \subseteq \mathbb{R}^2$ then 
$d_\mathcal{H}(A,B)<\epsilon$
for a small $\epsilon$ with average value obtained from our experimental data of $0.84\%$ of the peak to peak value of the signal. (The explicit construction process is explained in detail in the Appendix B.) On the other hand, for non-wheeze signals the Hausdorff distance between the graph of the model and that of the signal is very large, with an average of $9.7 \%$ of the peak to peak value of the signal according to our experimental dataset. This thus indicates that the proposed model is inadequate to represent normal signals. Comparison of the proposed model for one recorded non-wheeze signal and one recorded wheeze is presented in Fig. \ref{mdl}. 
%
%
\subsection{Topological Characterization of Delay Embeddings}
The continuous time delay embedding of $w(t)$ and $w_i(t)$ can be defined as
\begin{equation}
\begin{array}{ll}
W(t)=\lbrace \left(w(t),w(t+\tau)\right): t \in [0, \infty) \rbrace,\\
W_i(t)=\lbrace \left(w_i(t),w_i(t+\tau)\right): t_{i-1} \leq t, t+\tau \leq t_i \rbrace.
\end{array}
\end{equation}
In this subsection, we state some lemmas for the continuous time delay embedding of $w(t)$ with constant $A$ and then generalize the results to the case with time varying amplitude. We then sample  $W(t)$ to conclude that the first Betti number of the delay embedding of $w(t)$ obtained using (\ref{gde}) with the delay selected as discussed in Section \ref{delay} is always at least 1.
The following lemma shows that the delay embedding sets of two sinusoidal functions that differ by a phase are ``eventually equal''.
%
\begin{lem}\label{phase}
Suppose that $v_i(t)=A \sin\left(\frac{2 \pi}{T} t+\phi_i \right)$, $t \in [0,\infty)$, $i=1,2$, $\phi_2 > \phi_1$. Fix $\tau > 0$, and consider the delay embedding sets $V_i(t)=\lbrace \left( v_i(t),v_i(t+\tau) \right),t \in [0,\infty) \rbrace$, $i=1,2$.
Consider the parameter change $t'=t+\frac{\phi_2 -\phi_1}{2 \pi}T$ and construct the sets $V_i=\lbrace \left( v_i(t'),v_i(t'+\tau) \right),t' \in [\frac{\phi_2 -\phi_1}{2 \pi}T, \infty) \rbrace$ and $V_2=V_2(t)$. Then $V_1=V_2$.
%
\end{lem}
The phase term can be therefore ignored during the next analysis of delay embeddings.
The following lemma shows that for a careful choice of the delay $\tau$ the delay embedding set of a sinusoidal function is an ellipse.
The proofs of the following lemmas can be found in the Appendix.
\begin{lem}\label{ellipse}
Suppose that $u(t)=A \sin\left(\frac{2 \pi}{T} t+\phi\right), t \in [0,\infty)$ and $\tau\neq k \frac{T}{2}$ where $k\in \mathbb{Z}$. The continuous delay-coordinate embedding of $u(t)$ yields an ellipse with radii $\alpha, \beta =A \sqrt{1\pm\cos(\frac{2\pi}{T})\tau}$ 
centered at the origin with angle of rotation  $\pm45^{\circ}$. 
The length of the sides of the circumscribed square around that ellipse is equal to $2A$.
\end{lem}
Note that in lemma \ref{ellipse}, $T$ is fixed while $\tau$ can vary. 
The following lemma helps in analyzing the behavior of delay embeddings when the period varies in order to extend the results to a piecewise sinusoidal function with different periods.
\begin{lem}\label{chngvar}
The delay-coordinate embedding of $u_i(t)=A \sin(\frac{2\pi}{T_i} t)$ using delay $\tau_i$ can be expressed as
$U_i(t)=\left(A\sin\left(\frac{2\pi}{T_i}t\right),A \sin\left(\frac{2\pi}{T_i} (t+\tau_i)\right)\right), i=1,2$.
The delay-coordinate embedding sets $U_1(t)$ and $U_2(t)$ differ by a reparametrization if and only if 
$\frac{\tau_2}{\tau_1}=\frac{T_1}{T_2}.$
\end{lem}
\noindent Now, we define the discrete sampled finite version of the delay embeddings of $w(t)$ and $w_i(t)$ defined in (\ref{w})-(\ref{wi}) as follows
\begin{equation} \label{samp}
\begin{array}{ll}
W=\lbrace W(t_j): t_j \in I \rbrace \\
W_i=\lbrace W_i(t_j): t_j \in I,  t_j, t_j+\tau \in [t_{i-1}, t_i] \rbrace.
\end{array}
\end{equation}
If $K=\lbrace \left( w_i(t_j),w_{i+1}(t_j +\tau) \right): t_j \in I , t_j \in [t_{i-1} , t_i], t_j+\tau \in [t_i,t_{i+1}]\rbrace$ is the set of transition points coming from different $w_i$ segments, then
\begin{equation}
  W=\bigcup_{i=1}^n W_i\bigcup K
 \end{equation}
Note that the points in $K$ are close to $W_i's$.
This is a direct result of Equations (\ref{w})-(\ref{wi}),(\ref{samp}) and the definition of $K$.
%
%
%
%
 \begin{figure}[t]
      \centering
\includegraphics[width=1\linewidth]{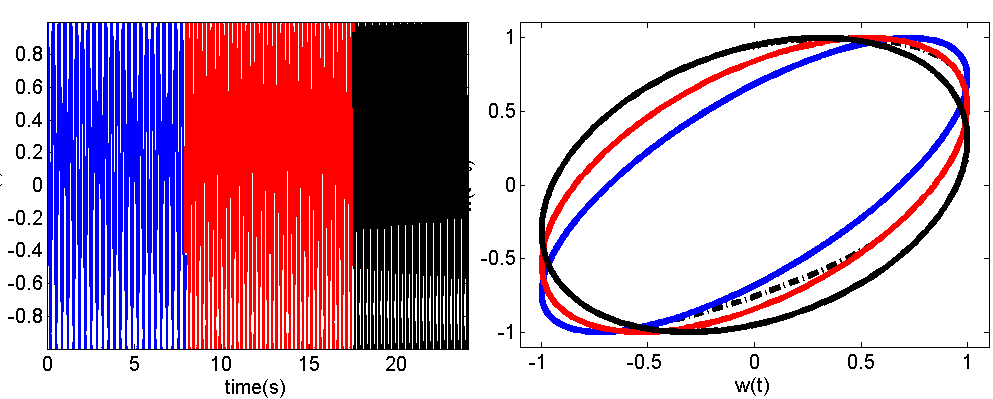} 
\caption{$w(t)$ with 3 different frequencies in different colors (left), the corresponding delay embedding $W(t)$ including the set $K$ shown using dashed curve (right). }
      \label{2d}
   \end{figure}
According to lemma \ref{ellipse}, for a fixed period, the radii of the ellipse depend only on the delay. Also, lemma \ref{chngvar} states that for a fixed period, changing the delay has the same effect on the delay embedding as a fixed delay and changing the period. Therefore, if we use a fixed delay and change the period, the radii of the ellipse are going to change accordingly. 
Thus, $\bigcup_{i=1}^n W_i$ 
with constant $A$ is a set of concentric ellipses with angle of rotation $\pm45^{\circ}$, and different radii inscribed in a unique circumscribed
square with side length of $2A$. This can be observed in Fig. \ref{2d}, where the delay embedding of a synthesized data with three different frequencies is illustrated. 
Note that
the finite set $K$ shown using dashed curve includes $1.6 \%$ of the whole data points.
%
Moreover, if the amplitude $A$ in (\ref{wi}) is time-varying, the delay-coordinate embedding $\bigcup_{i=1}^n W_i$  would be 
a set of concentric ellipses with an angle of rotation $\pm 45^{\circ}$ and different radii
inscribed in squares with different side lengths. 
\begin{lem}\label{le4}
Suppose that $\cup_{i=1}^n W_i$ is a union of concentric ellipses on the plane, with radii $\alpha_i,\beta_i=A_i\cdot\sqrt{1\pm\cos(\frac{2\pi}{T_i})\cdot\tau}$, centered around zero with angle of rotation $\pm 45^o$. Then
the persistence diagram of $\cup_{i=1}^n W_i$, has at least one 1-dimensional persistent bar.
\end{lem} 
Following the above discussion, one can conclude:
\begin{thm}\label{holemodel}
Suppose that the time delay $\tau$ is selected between $t_{c1}$ and $t_{c2}$. The delay-coordinate embedding $W$ is a set close to a set of concentric ellipses $\left(\bigcup_{i=1}^n W_i \right)$ with angles of rotation $\pm45^{\circ}$, different radii and different side lengths of the circumscribed squares around them. Therefore, the Rips complex associated to the corresponding point cloud of $W$ always has at least one 1-dimensional persistent hole.
\end{thm}
%
\noindent We will now utilize the obtained results to give a similar topological characterization for wheeze signals. For each digital wheeze signal $s(t)$, we construct a function $w(t)$ defined using (\ref{w})-(\ref{wi}) such that for $W$ as defined in (\ref{samp}) and $S$ defined in the same way for $s(t)$, we get
 $d_\mathcal{H}(S,W)<\delta$,
where $\delta=0.0119$ on average, according to our experimental results. 
%
%

The Wasserstein distance is a metric to measure the distance between the persistent barcodes. According to \cite{Ghzl}, the Wasserstein distance between persistent barcodes of two point clouds is less than or equal to the Hausdorff distance between the point clouds, i.e.
$W_{\infty}(\text{Bar}(S),\text{Bar}(W))\leq d_\mathcal{H}(S,W),$ 
where  $\text{Bar}(S)$ and $\text{Bar}(W)$ are the persistent barcodes of the finite point clouds $S$ and $W$, respectively. Thus,
\begin{equation}\label{Wdel}
W_{\infty}(\text{Bar}(S),\text{Bar}(W))< \delta. 
\end{equation}
Equation (\ref{Wdel}) shows that the persistent barcodes of the delay embedding point clouds of the wheeze signal and of the model are in $\delta-$distance of each other. Consequently,
%
%
%
the first Betti number of a point cloud of the delay embedding of a wheeze signal with the time delay selected as discussed in Section \ref{delay}
is always at least 1. On the other hand, using experimental results we show that the first persistent Betti number of the delay embedding of a non-wheeze signal is zero.
%
As a result of previous discussions, the use of this approach in the analysis techniques can address the problem of differentiating between wheeze and non-wheeze sound signals quite well.
\begin{figure}[t]
\centering
\begin{tabular}{cc}
\includegraphics[width=0.5\linewidth]{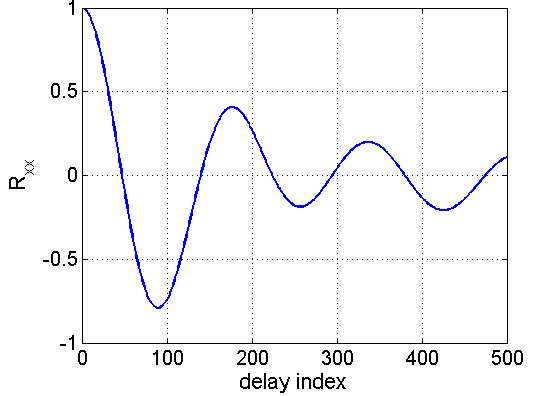} &
\includegraphics[width=0.5\linewidth]{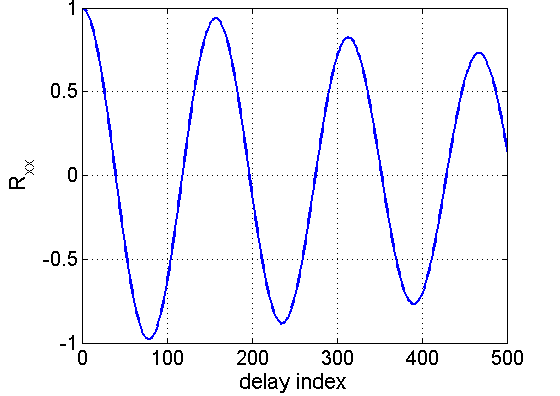}  \\
\end{tabular} \\
\caption{The autocorrelation like (ACL) function of a non-wheeze breathing sound signal (left) and a wheeze signal (right)} \label{cr}
\end{figure}
\section{Experimental Results}
\begin{figure}[t]
      \centering
      \includegraphics[width=0.95\linewidth]{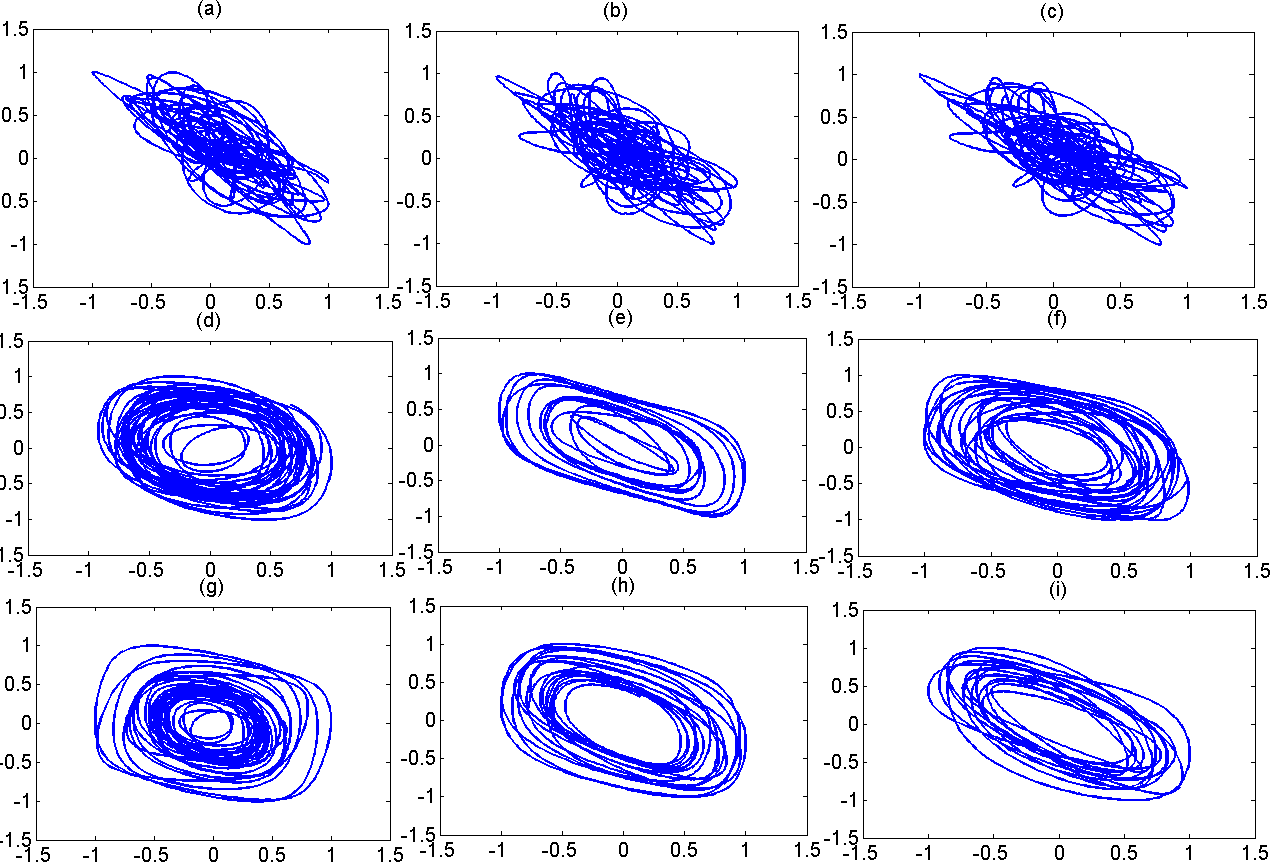}
      \caption{Delay-coordinate embedding of non-wheeze signals recorded over (a) apex, (b) midlung and (c) chest. Delay embeddings of different wheezes including (d) a high pitched (e) a fixed monophonic (f) a random monophonic (g) an expiratory monophonic wheeze (h) an expiratory, and (i) an inspiratory stridor}
      \label{dem}
   \end{figure}

We have performed experiments using a large number of recorded breathing sound signals taken from \cite{Lip, DW, RW} 
to support the effectiveness of the proposed approach. For presentation in this paper, we consider 9 different recorded breathing sound signals including 3 normal sounds recorded over different areas of the body
and 6 types of wheezes. The sampling rate is 44.1 kHz.
For each wheeze signal we choose a small interval where wheezes are heard by ear and confirmed by experts. We then normalize the amplitude between -1 and 1.
The zero of the ACL function (Fig. \ref{cr}) of the signal is selected as a delay in order to obtain the maximum size hole inside the point cloud if there exists any. Note that since breathing sound signals are time varying and non-stationary, a strict autocorrelation function cannot be used for periodicity analysis.
Fig. \ref{dem} represents the 2D time delay-coordinate embeddings with the chosen delay parameters. 
  \begin{figure}[t]
      \centering
      \includegraphics[width=0.95\linewidth]{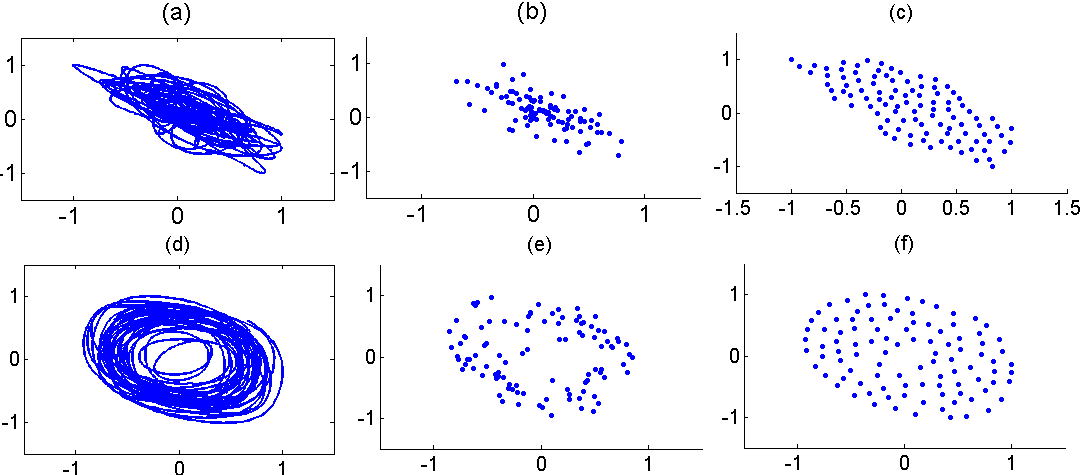}
      \caption{Point cloud subsampling, (a) and (d): the delay embedding of a non-wheeze and a wheeze signal including 4000 points. 100 subsamples selected using random method ((b),(e)) and maxmin method((c),(f))}
      \label{ln}
   \end{figure}
\begin{figure}[tb]
\centering
\includegraphics[width=0.95\linewidth]{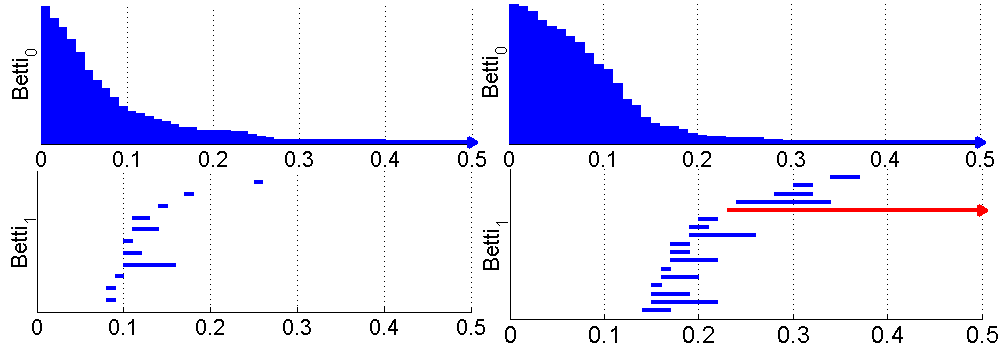}
\caption{The barcodes for a non-wheeze breathing sound signal (left) and a wheeze signal (right). The ``significant" barcode is highlighted in red and is used to distinguish wheeze signals from non-wheeze signals.} \label{bc}
\end{figure}
For each point cloud, a Rips complex is constructed using javaPlex \cite{javaPlex} as follows: We first select 100 subsamples using  maxmin and random method described in section \ref{landmark selection} which are about $3\%$ of the whole number of points. Fig. \ref{ln} shows subsamples selected from the point cloud of one non-wheeze and one wheeze signal shown in Fig. \ref{dem}-(a) and (d), respectively. 
It is easier to detect the hole when using random subsamples for its lower computational complexity.
The computation of $ \beta_0$ and $\beta_1$ is carried out and shown in Fig. \ref{bc} by the barcodes 
for the random subsamples of the non-wheeze and wheeze signals as shown in Fig. \ref{ln}. 
The accuracy of our proposed technique 
is $98.39\%$ while the accuracy of the techniques proposed in \cite{Homs04} and \cite {Taplidou07} are  $86.2\%$ and $95.5\%$, respectively.
\section{Conclusion}
We proposed a persistent homology analysis of a delay coordinate embedding of a point cloud to capture almost harmonic behaviours.
%
A very small number of data points are selected from all points to represent the sound signal. Using a subsampling method reduces the computational complexity of the algorithm making it implementable in low-power systems. 
\section{Appendix}
\subsection{Mathematical Proofs}

\begin{proof}
(Lemma 2) We first assume that $A=1$ for simplicity, then generalize the proof for any constant value of $A$. The trigonometric parametrization of the ellipse can be used for proof. The general quadratic equation $ax^2+bxy+cy^2-\rho=0$ represents an ellipse with center at the origin if the following conditions are satisfied:
\begin{equation}\label{23}
\left\lbrace \begin{array}{ll}
  b^2-4ac<0 \\
   \delta=\begin{vmatrix} 
2a & b & 0 \\
b & 2a & 0\\
0 & 0 & -2\rho
\end{vmatrix}\neq0,
\\
       \delta.\, \sigma < 0, \quad \text{where}\,\,\,\sigma= a+c.
   \end{array}
   \right.                             
\end{equation}
%
\noindent Moreover, the angle of rotation of the ellipse can be obtained using the following equation
\begin{equation}\label{angle}
\tan(2\theta)=\frac{b}{a-c}.
\end{equation}
For notational convenience let $C:=\cos\theta$ and $S:=\sin\theta$. The radii of the ellipse are then
\begin{equation}\label{radii}
\alpha=\sqrt{\frac{\rho}{P}},\quad   \beta=\sqrt{\frac{\rho}{Q}},
\end{equation}
where $P=aC^2+cS^2+bCS$ and $Q=aS^2+cC^2-bCS$.
According to lemma 3, the delay coordinate embedding of $A\sin(\frac{2\pi}{T} t)$ with delay $\tau$ is equivalent to the delay embedding of $A\sin(t)$ with the delay $ \frac{2\pi}{T} \tau$. Therefore, if we prove that the delay-coordinate embedding of $\sin(t)$ constructs an ellipse centered at the origin with the angle of rotation  $45^{\circ}$ provided that $\tau\neq k \pi$ where $k\in \mathbb{Z}$, the first part of the lemma will be proved. 
Let us construct the quadratic equation using $x(t)=\sin(t)$ and its delayed version $y(t)=\sin(t+\tau)$. Simplifying the obtained equation yields,
\begin{equation}\label{quad}
\begin{array}{l}
{\sin^2t}\left(a+b \cos \tau + c \,{\cos^2 \tau}\right) + \cos^2 t\left( c\, \sin^2 \tau \right)  \\ 
\qquad \qquad   \qquad   +\cos t \sin t\left( b\, \sin \tau + 2c\, \sin \tau \cos \tau \right)- \rho  = 0.
\end{array}
\end{equation} 
In order for this equation to be true for all $t$, the following conditions should be satisfied
\begin{equation}\label{22}
\left\lbrace  \begin{array}{ll}
                                 a + b\cos\tau+ c \cos^2\tau =c \sin^2\tau,\\
                              b \sin\tau +2 c \sin\tau \cos\tau=0.
                             \end{array}
                             \right.                              
\end{equation}
Simplifying Equation (\ref{22}) yields 
\begin{equation}\label{23}
\left\lbrace  \begin{array}{ll}
                                 a = c,\\
                               b + 2 c \cos\tau=0.
                             \end{array}
                             \right.                             
\end{equation}
In this case, Equation (\ref{quad}) becomes
\begin{equation}\label{rho}
c \sin^2\tau-\rho=0.
\end{equation}
Thus we can check the conditions. We can first calculate
\begin{equation}
b^2-4ac=(-2 c \cos\tau)^2-4c^2=4c^2\left(\cos^2\tau-1\right)
\end{equation}
First conditions are therefore satisfied since $\tau\neq k\pi$,  $k\in \mathbb{Z}$. 
According to (\ref{23}) and (\ref{rho}), one can calculate $\delta$ as
\begin{equation}
\begin{array}{l}
\delta=-2\rho\left(4a^2-b^2\right)=-8\rho a^2(1-\cos^2\tau)\\
\qquad \qquad   \qquad \qquad \, \,=-8a^3\left(1-\cos^2 \tau\right)^2,
\end{array}
\end{equation}
which can never be zero since $\tau\neq k\pi$, $k\in \mathbb{Z}$. Moreover, it is of opposite sign to $a$ and thus of opposite sign of $\sigma$, since $\sigma=2a$. Hence, all three conditions are satisfied and the delay coordinate embedding of $\sin(t)$ with delay $\tau\neq k\pi$,$\,k\in \mathbb{Z}$ defines an ellipse centered at the origin. Additionally, 
according to (\ref{angle}), the angle of rotation is $45^{\circ}$. The same results are also true with delay $\frac{2\pi}{T} \tau$, where $\tau\neq \frac{kT}{2}$.
Therefore, $C=S=\sqrt{2}/2$.
Now, the radii of the ellipse can be calculated since
\begin{equation}
\begin{array}{cc}
                                 P=(a+b+c)/2=a(1-\cos\tau),\\
                               Q=(a-b+c)/2=a(1+\cos\tau),
                             \end{array}
\end{equation}
$\alpha=\sqrt{1+\cos\tau}$ and $\beta=\sqrt{1-\cos\tau}$.
\noindent Clearly, $\alpha,\beta\neq0$ since $\tau\neq k\pi$. So if we change the delay to $ \frac{2\pi}{T} \tau$, the radii of the ellipse are going to be $\alpha=\sqrt{1+\cos(\frac{2\pi}{T}) \tau}$ and $\beta=\sqrt{1-\cos(\frac{2\pi}{T} \tau)}$. Multiplying by the amplitude does not change the conditions since it just affects $\rho$, and multiplies it by $A^2$ in (\ref{quad}). Therefore the radii of the ellipse would be multiplied by $A$ according to (\ref{radii}).
Since the delay coordinate embedding of $A\sin(\frac{2\pi}{T} t)$ with delay $\tau$ is equivalent to the delay embedding of $A \sin(t)$ with delay $\frac{2\pi}{T} \tau$, the lemma is proved.
\end{proof}

\begin{proof}
(Lemma 3) Using a change of variable $t'=\frac{T_1}{T_2}t$, 
\begin{equation}
\begin{aligned}
U_2(t')=\left( A \sin\left(\frac{2\pi}{T_1} t'\right),A \sin\left( \frac{2\pi}{T_1}t'+\frac{2\pi}{T_2}\tau_1\right)\right)\\
=\left(A \sin \left( \frac{2\pi}{T_1} t'\right),A \sin\left( \frac{2\pi}{T_1}t'+\frac{T_1}{T_2}\tau_1\right)\right)
\end{aligned}
\end{equation}
So, $U_1(t)=U_2(t')$ for $ \tau_2=\frac{T_1}{T_2}\tau_1$. Thus, the sets $U_1(t)$ and $U_2(t)$ contain the same points implying that the two sets are equal.
\end{proof}

\begin{proof}
(Lemma 4)
Let $W_s$ be the ellipse with the smallest radii
\begin{equation}
\alpha_s,\beta_s = A_s\sqrt{1\pm\cos(\frac{2\pi}{T_s})\tau}.
\end{equation}
In the continuous case, the persistent diagram for $\cup_{i=1}^n W_i$ would contain a persistent bar with birth time $t_b=0$ and death time $t_d=\beta_s$ corresponding to this small ellipse. In the discrete case, the birth time of this bar will be close to $t_b=0$ depending on the $T_s$, and the death time $t_d \geq \beta_s$. Therefore, the persistent diagram has a bar of length close to $\beta_s$.
\end{proof}

\begin{proof}
(Theorem 1) The delay $\tau$ is chosen between $t_{c1}$ and $t_{c2}$. According to our experimental results, selecting the second zero of ACL function gives the best informative delay embedding. Similar to the proof of lemma 4, the smallest ellipse $W_s$ will create a significant 1-dimensional persistent bar of length close to $\beta_s$. The points in $K$ do not interfere with this ellipse since they are in between $W_i's$ and close to them (Fig. 2). 
\end{proof}


\subsection{The model construction}\label{const}
\noindent We construct a model corresponding to each wheeze signal by generating a signal with similar frequency and amplitude as those of the wheeze signal at each time instance. This objective is achieved in two steps. First, we estimate the frequency of the signals at different time intervals using the zero crossing points. The intervals between successive zero crossing indexes can be considered as a random variable $d_i, i=1,2,...,n_z-1$ where $n_z$ is the number of zero crossings .
We can find the time intervals at which the frequency is constant using the distribution of $d_i$, and detect the intervals with almost constant frequencies denoted by $T_j,j=1,2,..,n_f$, where $n_f$ is the number of different frequencies of the signal. We then find the frequency during each time interval as  $f_j=\sfrac{1}{2\mu_j}$, where $ \mu_j,j=1,2,..n_f$ denotes the mean value of $d_i$ during $T_j$.
Next, we find the envelope of the signal using its critical points and multiply it by a sinusoidal function with the frequencies obtained in the first step. Only critical points with positive amplitude are used since the breathing sound signals are to a large extent, almost symmetric about the time axis.
Next, we perform an interpolation between them to obtain the envelope of the signal denoted by $a(t)$ and multiply it by the obtained piecewise sinusoidal function.

\bibliographystyle{ieeetr}
\bibliography{Arxive4}
\end{document}